%% file: DarioTom12DB.tex
\documentclass{article}
\usepackage{amsmath} 
\usepackage{amsfonts,mathrsfs}
\usepackage{amssymb}
\usepackage{color}
\usepackage{graphicx}
\usepackage{epsfig}
\usepackage{subfigure}

\newtheorem{definition}{Definition}{\bf}{}
{\bf}{}
\newtheorem{theorem}{Theorem}{\bf}{}
\newtheorem{prop}{Proposition}{\bf}{}
{\bf}{}
\newtheorem{assumption}{Assumption}{\bf}{}
\newcommand{\proof}{\vspace{1mm}\noindent{\it Proof.}\quad}
 
\newtheorem{problem}{Problem}
\usepackage{color}
\newtheorem{example}{Example}

\newcommand{\dR}{\mathbb{R}}

\newcommand{\rmd}{{\rm{d}}}

\renewcommand{\phi }{\varphi }
\def\dist{{\rm dist}}

\title{Approachability in Population Games\thanks{The work of D.\ Bauso was supported by the 2012 ``Research Fellow'' Program of the Dipartimento di Matematica, Universit\`a di Trento and by PRIN 20103S5RN3 ``Robust decision making in markets and organizations, 2013-2016''.}}
\author{Dario Bauso\thanks{D.\ Bauso is with Dipartimento di Ingegneria Chimica, Gestionale, Informatica e Meccanica, Universit\`a di Palermo, Italy. D.\ Bauso is currently visiting professor at the Department of Engineering Science, University of Oxford, UK, email: \textsl{dario.bauso@unipa.it}} \and Thomas W.\ L.\ Norman\thanks{T.\ Norman is at Magdalen College, Oxford, email: \textsl{thomas.norman@magd.ox.ac.uk}}}

\begin{document}
\maketitle

\begin{abstract}
This paper reframes approachability theory within the context of population games. Thus, whilst one player aims at driving her average payoff to a predefined set, her opponent is not malevolent but rather extracted randomly from a population of individuals with given distribution on actions. 
First, convergence conditions are revisited based on the common prior on the population distribution, and we define the notion of \emph{1st-moment approachability}. Second, we develop a model of two coupled partial differential equations (PDEs) in the spirit of mean-field game theory:\ one describing the best-response  of every player given the population distribution (this is a \emph{Hamilton-Jacobi-Bellman equation}), the other capturing the macroscopic evolution of average payoffs if every player plays its best response (this is an \emph{advection equation}). Third, we provide a detailed analysis of existence, nonuniqueness, and stability of equilibria (fixed points of the two PDEs). Fourth, we apply the model to regret-based dynamics, and use it to establish convergence to Bayesian equilibrium under incomplete information.
\end{abstract}
\section{Introduction}\label{sec:introduction}

We consider a game played by a large population of individuals in continuous time. At every time, each individual engages in play with a random opponent  extracted from the population and  the resulting payoff, which depends on the action profiles of both players, is a vector. Such vector payoffs can be interpreted as deriving from a collection of noninterchangeable goods. Let us think, for instance, of a negotiation between an employer and a candidate employee over salary, career prospects, maximal number of days off and so forth. Formally, we can think of the completeness axiom being satisfied along each dimension of our vector but failing across them, giving a special case of Aumann's \cite{A62} 
framework. Indeed, vector payoffs may also be appropriate when the continuity axiom fails (see \cite{BBD91}). 
Alternatively, each player may be representative of a group of individuals whose preferences may not be aggregated into a single ordering, so that the vector payoff has one component for each individual in the group
. Finally, payoff vectors also naturally arise when considering their regret at not having made each possible deviation.

\noindent
\textbf{Main results.} First, we provide a new model that combines approachability and population games. 
Given that the opponent is randomly extracted from the population, the approach by Blackwell---which looks at the worst-case payoff---may appear conservative. Thus, we relax Blackwell's conditions, assuming that the opponent is not malevolent but instead is simply extracted from a population with given distribution; we call this \emph{$1$st-moment approachability}. Second, we build upon the theory of mean-field games and adapt the concept of mean-field equilibrium to our evolutionary set-up; we call this \emph{self-confirmed equilibrium}. Third, we discuss existence and nonuniqueness of the equilibrium. Finally, we explore the regret interpretation of our model; whereas $1$st-moment approachability of nonpositive regrets no longer implies Nash equilibrium (as in \cite{HM03}), we show that nonpositive maximal regret does imply Bayesian equilibrium under incomplete information.

\noindent
\textbf{Related literature.} The theory of ``approachability'' dates back to Blackwell \cite{B56} and culminates in the well known Blackwell's Theorem. Approachability arises in several areas of game theory, such as allocation processes in coalitional games \cite{L02-CT}, regret minimization \cite{L03,HM03},
adaptive learning \cite{CL06,FV99,H05,HM01}, excludability and bounded recall \cite{LS06}, and weak approachability \cite{V92}, just to name a few.
For instance, in coalitional games one asks whether the core is an approachable set, and which allocation processes can drive  the complaint vector to that set.  In regret minimization, one considers the nonpositive orthant in the space of regrets; a player tries to adjust her strategy based on the current regret so as to make that set approachable by the regret vector. Once all players have nonpositive regret, the resulting outcome is an equilibrium for the game. This idea of adapting the new action to the current state of the game is common to adaptive learning and evolutionary games as well, but in regret-based dynamics the state is in payoff (rather than strategy) space. Evolution under incomplete information has been relatively little studied, with the notable exception of Ely and Sandholm \cite{ES05,S07}, who analyse a best response dynamic with a subpopulation for each possible type; here, by contrast, we have a single population of agents with nonconstant types who adopt (type-dependent) Bayesian strategies through time.

Despite its discrete-time nature in the original Blackwell formulation, approachability has been extended to continuous-time repeated games, thus showing common elements with Lyapunov theory \cite{HM03}. Though first formalized in finite-dimensional spaces, a definition of approachability in infinite-dimensional space has been provided by Lehrer \cite{L02}.
Approachability can be reframed within differential games and as such can be studied using differential calculus and stability theory  \cite{LS07,SQS09}. In  particular, in \cite{LS07} the authors show that, beyond being an extension (to a vector space) of the von Neumann minmax theorem \cite{vN28}, the approachability principle also has elements in common with differential inclusions \cite{AC84}.
In addition to this, \cite{SQS09} establishes connections with viability theory \cite{A91}, and set-valued analysis \cite{AF90} (see, cfg., the comparison between an approachable set and  a discriminating set) and set invariance theory \cite{B99}.\footnote{Still within the realm of differential games, it is worth noting that the notion of nonanticipative behavior strategies has a long history \cite{BB11,EK72,SQS09,R69,V67}. Actually, it turns out that classical feedback strategies in differential games are special nonanticipative strategies.}

The approachability principle is also behind the notion of excludability; along this line, some authors investigate which sets are approachable and which ones excludable under imperfect information (bounded recall, delayed and/or stochastic monitoring) \cite{LS06}. 
Connected to approachability as well is the concept of ``attainability.'' 
Attainability is a new notion developed in \cite{BLS12,LSB11} in the context of 2-player continuous-time repeated games with vector payoffs. Attainability arises in many contexts such as transportation networks, distribution networks, production networks applications. The main question is: ``Under what conditions does a strategy for player 1 exist such that the cumulative payoff converges (in the lim sup sense) to a preassigned set (in the space of vector payoffs) independently of the strategy used by player 2?''

A second stream of literature we follow in the present study is the one on \emph{mean field games}. This theory originated in the work of M.\ Y.\ Huang, P.\ E.\ Caines and R.\ Malham\'e \cite{HCM03,HCM06,HCM07}, and independently in that of J.\ M.\ Lasry and P.\ L.\ Lions \cite{LL06a,LL06b,LL07}, where the now standard terminology of Mean Field Games (MFG) was introduced.
 Explicit solutions in terms of mean field equilibria are not common unless the problem has a linear-quadratic structure, see \cite{B12}.
 Mean field games have connections to evolutionary games (see for instance \cite{javano}) and large games \cite{A64}. Actually, both the \emph{anonymous game} in \cite{javano} and the \emph{large game} in \cite{A64} build upon the notion of mass interaction and can be seen as a stationary mean field.

This paper is organized as follows. 
In Section \ref{sec:setup}, we set up the problem. In Section \ref{motivations},
we provide our population game motivation for the problem at hand.   
In Section \ref{main}, we establish the main results of the paper. 
In Section \ref{regret}, we apply the model to a regret-based setting, and show under incomplete information that nonpositive maximal regrets that are approachable in $1$st moment must be Bayesian equilibria. Finally, in Section \ref{sec:conclusions}, we draw  concluding remarks. 

\noindent
{\bf Notation}. We view vectors as columns.  For a vector $x$, we use $x_i$  to denote
its $i$th coordinate component. Occasionally we may write $(x)_{i=1,\ldots,m}$ to denote an $m$-dimensional column vector. For two vectors $x$ and $y$, we use $x<y$ ($x\le y$) to denote
$x_i<y_i$ ($x_i\le y_i$) for all coordinate indices $i$. 
We let $x^T$ denote the transpose of a vector $x$, and $\|x\|$ its Euclidean norm.  We write $P(x)$ to denote the projection of a vector $x$ on a set $X$,
and $\dist(x,X)$ for the distance from $x$ to $X$, i.e.\ 
$P(x) = \arg \min_{y \in X} \|x - y\|$ and $\dist(x,X)=\|x-P(x)\|$,
respectively. We also denote by $conv$ the convex hull of a given set of points.
$\partial_x$ indicates the first partial derivative with respect to $x$.

\section{The Model}\label{sec:setup}
With the above preamble in mind, the game at hand is a two-player repeated game with vector payoffs in continuous time.\footnote{Whilst the game is repeated, opponents are constantly rematched, and hence no supergame considerations arise.} We assume that the players use nonanticipating behavior strategies with delay. This means that the behavior of a player may depend only on past play. In other words, the way a player plays during a given interval of time does
not affect the way the  opponent  plays during that block. Still,
it may affect the other player's play in subsequent intervals.

Let $A=\{1,2,\ldots,n\}$ be a discrete set, $a_i:[0,T] \to A$ a measurable function of time and $a_j:[0,T] \to A$ a random disturbance. 
 Let  $u: A \times A \to M$ where $M=\{M_{lk}, l,k\in A\}$   and   $M_{lk} \in \dR^m$ (each entry $M_{lk}$ is an $m$-dimensional vector). 
Let $X:=conv\{M_{lk}|\, l,k \in A\}$, where $conv$ denotes the \emph{convex hull}, and  consider the 
differential equation in $X$
\begin{equation}\label{dyng}
\left\{\begin{array}{ll}
d x(t) = \frac{1}{t} (\mathbb E u(a_i(t),a_j(t)) - x(t)) dt,\quad \forall t \in[0,T],\\ 
x(0)=x_{0}\in X,
\end{array}\right.
\end{equation}
where $x_{0}$ is generated according to a distribution law $\rho_0(x)$. 
More specifically, consider a  probability density function $\rho: X \times [0,+\infty[ \to \mathbb R$, $(x,t) \mapsto \rho(x,t)$, representing the density of the players whose state is $x$ at time $t$, which satisfies $\int_{\mathbb R} \rho(x,t) dx=1$ for every $t$.  Let us also define the mean state over players at time $t$  as $\overline \rho(t) := \int_X x \rho(x,t) dx$. We also have $\rho(x,0)=\rho_0(x)$.

The objective of a player  is to approach a given target $y:[0,T] \to X$. Then, for each group, consider a  running  cost  $g:X \times X \to [0,+\infty[$, $(x, y)\mapsto g(x,y)$  of the form: 
\begin{eqnarray}\label{gg} 
g(x, y) & = & \frac{1}{2}\left[\left(y - x\right)^T Q \left(y - x\right)
 \right],
\end{eqnarray}
where $Q>0$ 
 and symmetric.

The above cost describes i) the  (weighted) square deviation of an individual's state from the target.

Also consider a terminal cost $\Psi:X \times X \to[0,+\infty[$, $(x,y)\mapsto\Psi(x,y)$ of the form
\begin{equation}\label{psig}\Psi(x,y) = \frac{1}{2} (y - x)^TS(y - x),\end{equation}
where $S>0$.
The problem in its generic form is then the following:

\begin{problem}\label{prob1}


Let the initial state $x(0)$ be 
 given and with density $\rho_0.$  
Given a finite horizon $T>0$, 
 a suitable running cost: $g: X \times X  \to [0,+\infty[$, $(x,y)\mapsto g(x,y)$, as in (\ref{gg}); a terminal cost $\Psi: X \times X \to [0,+\infty[$, $(y,x)\mapsto\Psi(y,x)$, as in (\ref{psig}), and  given a suitable dynamics for $x$ as in (\ref{dyng}), solve 
\begin{eqnarray}
\inf_{a_i(\cdot) \in \mathcal C}  \left\{J(x_0,a_i(\cdot),a_j(\cdot))
=  \int_0^T g(x(t),y)dt + \Psi(x(T),y)\right\},\label{cost}
\end{eqnarray} 
\noindent
where  $\mathcal C$ is the set of all measurable functions $a_i(\cdot)$ from $[0,+\infty[$ to $A_i$, and $\mathbb E u(\cdot)$ in (\ref{dyng}) must be consistent with the evolution of the distribution $\rho(\cdot)$ if every player behaves optimally. 
\end{problem}

\section{Motivation: Population Games}\label{motivations}

Consider a population game where continuously in time every individual matches with an opponent randomly extracted from the population and the resulting payoff is a vector. 
The resulting game is a two-player repeated game with vector payoffs in continuous time $\Gamma$ that every individual plays against a population with given (evolving) distribution over actions. 
Let $A$ be the finite set of actions of every individual, then the instantaneous payoff is given by a function $u: A \times A \to \dR^m$, where $m $ is a natural number. We assume w.l.o.g.\ that payoffs are bounded and correspond to the elements of the following discrete set $M=\{M_{lk}, l,k\in A\}$ where $M_{lk} \in \dR^m$, so that $u: A\times A \to M$. We extend $u$ to the set of mixed-action pairs, $\Delta(A) \times \Delta(A)$, in a bilinear fashion. The one-shot vector-payoff game $(A,A,u)$ is
denoted by $G$ and we will say that the game in continuous time $\Gamma$ is \emph{based on} $G$.

The game $\Gamma$ is played over the time interval $[0,\infty)$.
We assume that the players use markovian strategies 
$$\sigma: X \times [0,T] \to A \, \mbox{  such that }\, a_i(t):=\sigma(x,t),$$ where $X:=conv\{M_{lk}|\, l,k \in A\}$ and $x$ is the average (over time) expected (over opponent's play) payoff defined as:
\begin{equation}
\label{equ966} 
x(t) = \frac{1}{t}\int_0^t \mathbb E u(a_i(t),a_j(t)) \rmd t \in
\dR^m
\end{equation}
In the above equation, 
\begin{equation}\label{expval0}
\left\{
\begin{array}{lll}
\mathbb E u(a_i(t), a_j(t))& := &u(a_i(t), q(t))\\
\\
q \in \Delta(A) \, & s.t. & q_k = \int_{R_k} \rho(x,t) dx, \\ &&  R_k:= \{x \in \mathbb R^m| \, \sigma(x,t)=k\}, \, \forall k \in A.\end{array}\right.
\end{equation}

Once we differentiate (\ref{equ966}) with respect to $t$ we obtain the equation 
(\ref{dyng}) in the same spirit as in Hart and Mas-Colell's paper \cite{HM03} on continuous-time approachability.
Then, Problem \ref{prob1} analyzes the approachability of a given target in the space of vector payoffs on the part of a population of individuals. 
 
\begin{figure} [htb]
\centering
\def\svgwidth{0.6\columnwidth}
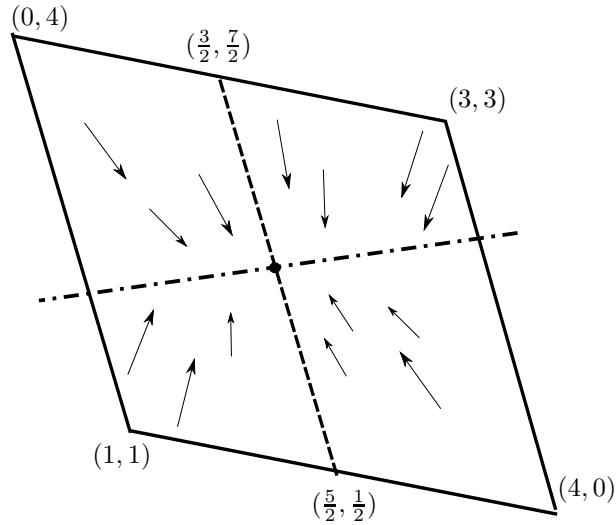
\caption{Payoff space of Prisoners' dilemma:\ State space $X=conv\{(3,3),(1,1),(0,4),(4,0)\}$ (boundary in solid line), supporting hyperplane $H$ (dot-dashed line) passing through the barycenter, vector field $dx(t)$ converging to $(\frac{3}{2},\frac{7}{2})$ for those who cooperate (region below $H$) and to $(\frac{5}{2},\frac{1}{2})$ for those who defect (region above $H$), $conv\{(\frac{3}{2},\frac{7}{2}),(\frac{5}{2},\frac{1}{2})\}$ is set of approachable points with population strategy $q=((\frac{1}{2},\frac{1}{2}),(\frac{1}{2},\frac{1}{2}))$, barycenter is self-confirmed with uniform distribution over $X$.}  \label{fig:PDp}
\end{figure}

\begin{example} \textbf{(Prisoners' Dilemma)}
Suppose, for instance, that players target the average payoffs across the population. Consider the following game:
$$\begin{array}{cc}
\begin{tabular}{|c|c|c|}
  \hline
       & Cooperate & Defect \\
  \hline
  Cooperate & $(3,3)$ & $(0,4)$ \\
  \hline
  Defect & $(4,0)$ & $(1,1)$\\
  \hline
\end{tabular}
\end{array}$$
Figure \ref{fig:PDp} depicts the payoff space in the continuous-time game based on this Prisoner's Dilemma. Here, the state space is  $X=conv\{(3,3),(1,1),(0,4),(4,0)\}$ (the boundary is in solid line), and the target $y$ is the barycenter assuming a uniform distribution.  One can visualize the supporting hyperplane $H$ (dot-dashed line) passing through the barycenter, and the vector field $dx(t)$ converging to $(\frac{3}{2},\frac{7}{2})$ for those who cooperate (region below $H$) and to $(\frac{5}{2},\frac{1}{2})$ for those who defect (region above $H$). The set $conv\{(\frac{3}{2},\frac{7}{2}),(\frac{5}{2},\frac{1}{2})\}$ is the set of approachable points with population strategy $q=((\frac{1}{2},\frac{1}{2}),(\frac{1}{2},\frac{1}{2}))$, and the barycenter  is at the equilibrium with uniform distribution over $X$.  This will be explained in Theorem \ref{polys}.
\end{example}

\section{Main results}\label{main}
This section outlines the main results of this paper. After introducing the \emph{expected value of the projected game}, Theorem \ref{app_princ} establishes conditions for approachability in 1st-moment. Theorem \ref{polys} introduces the notion of self-confirmed equilibrium. Theorems \ref{thm:ex} and \ref{thm:uq} elaborate on existence and nonuniqueness respectively.

\subsection{Expected value of the projected game}
Given the above game, we wish to analyze convergence properties in the space of distributions of the cumulative or average payoff $x_i(t)$,  in the spirit of  approachability. We will make use of the notion of \emph{projected game} which we recall next. 
Let $\lambda \in \mathbb R^m$ and denote by $\langle \lambda,G\rangle$ the one-shot zero sum game whose set of players and their actions are as in game $G$, and the payoff that player $j$ pays to player $i$ is $\lambda^T u(a_i(t),a_j(t))$ for every $(a_i(t),a_j(t)) \in A_i \times A_j$. Observe that, as a zero-sum one-shot game, the game $\langle \lambda, G\rangle$ has a \emph{value},  $val(\lambda)$, obtained as $$val(\lambda):=\min_{a_i(t)} \max_{a_j(t)} \lambda^T u(a_i(t),a_j(t)).$$ 
Given the stochastic nature of $a_j(t)$ the above min-max operation is not useful to our purposes. Then, we  rather consider the expected value of the game (where the inner maximization is replaced by an expectation) 
 and discuss approachability in expectation. 
 In the light of this, and using the bilinear structure of the utility function, and assuming markovian strategies $$\sigma: X \times [0,T] \to A \, \mbox{  such that }\, a_i(t):=\sigma(x,t)$$ we can rewrite the expected value as
\begin{equation}\label{expval}
\left\{
\begin{array}{lll}
\mathbb E val(\lambda)& := & \min_{a_i(t)} \mathbb E \lambda^T u(a_i(t), a_j(t))\\ &=&\min_{a_i(t)}  \lambda^T u (a_i(t),  q(t) ),\\
\\
q \in \Delta(A) \, & s.t. & q_k = \int_{R_k} \rho(x,t) dx, \\ &&  R_k:= \{x \in \mathbb R^m| \, \sigma(x,t)=k\}, \, \forall k \in A.\end{array}\right.
\end{equation}
In the case of state-dependent payoff, which occurs when we consider the game whose payoff is $$f(u(a_i(t),a_j(t)), x(t)) = \frac{1}{t} (\mathbb E u(a_i(t),a_j(t)) - x(t) )= \frac{1}{t} (u(a_i(t),q(t)) - x(t) ),$$  the above expression can be  modified as:
\begin{equation}\label{expvalsd}
\left\{\begin{array}{lll}
\mathbb E val_x(\lambda)& := & \min_{a_i(t)} \mathbb E \lambda^T f\Big(u(a_i(t), a_j(t)),x_i)\Big)\\ &=&\min_{a_i(t)}  \lambda^T f\Big(u (a_i(t),q(t) ),x_i \Big)\\
\\
q \in \Delta(A) \, & s.t. & q_k = \int_{R_k} \rho(x,t) dx, \\ &&  R_k:= \{x \in \mathbb R^m| \, \sigma(x,t)=k\}, \, \forall k \in A.\end{array}\right.
\end{equation}
Note that here we use the notation $u(a_i(t),q(t))$ to mean $\mathbb E u(a_i(t),a_j(t))$.

\subsection{Approachability in 1st-moment}
Approachability theory was developed by Blackwell in 1956 \cite{B56}
and is captured in the well known Blackwell's Theorem. 
We recall next the geometric (approachability) principle that lies behind Blackwell's Theorem. 

To introduce the approachability principle, let $\Phi$ be a closed and convex set in $\mathbb R^m$ and let $P(x)$ be the projection of any point $x \in \mathbb R^m$ (closest point to $x$ in $\Phi$).

\begin{definition}(Approachable set) 
A closed and convex set $\Phi$  in $\mathbb R^m$ is \emph{approachable} by player 1 if there exists a strategy for player 1 such that (\ref{appr11}) holds true for every strategy  of player 2:
\begin{equation}\label{appr11}\lim_{t \rightarrow \infty} dist(x(t),\Phi)=0.\end{equation}\end{definition}
The next result is the Blackwell approachability theorem.
\begin{prop}\label{app_princ}(Approachability principle \cite{B56,LS07}) A closed and convex set $\Phi$  in $\mathbb R^m$ is approachable by player 1 if for every
$x(t)$ there exists a strategy for player 1 such that (\ref{appr}) holds true for every strategy of player 2:
\begin{equation}\label{appr} [x(t) - P(x(t))]^T [x(t)-P(x(t)) + f(u_i(\sigma(x,t),a_j(t)), x_i(t))] \leq 0, \quad \forall \ t. \end{equation}
\end{prop}

Note that in the above statement, condition (\ref{appr}) is equivalent to saying that i) for every $x$ 
taking $\lambda = \frac{x-P(x)}{\|x-P(x)\|} \in \mathbb R^m$ the value of the projected game satisfies

\begin{equation}\label{appr1} [x(t) - P(x(t))]^T [x(t)-P(x(t))] + \|x-P(x)\| val_x(\lambda) \leq 0, \quad \forall \ t. \end{equation}

Now, if we assume that the opponent is committed to play a mixed strategy $q \in \Delta(A)$,
 condition (\ref{appr}) turns into 
\begin{equation}\label{apprm}[x(t) - P(x(t))]^T [x(t)-P(x(t)) + f(u(\sigma(x,t),q(t)), x(t))]\leq 0, \quad \forall \ t,\end{equation}
and the corresponding condition (\ref{appr1}) can be rewritten as 
\begin{equation}
\left\{\begin{array}{lll}[x(t) - P(x(t))]^T [x(t)-P(x(t))]+ \|x-P(x)\| \mathbb E val_x(\lambda)\leq 0, \quad \forall \ t,\\
\mathbb E val_x(\lambda) := 
 \min_{a_i(t)}  \lambda^T f (u_i (a_i(t),q(t) ),x_i).\end{array}\right.
\end{equation}

\begin{theorem}\label{app_princ} \textbf{(Approachability in $1$st-moment)}
Let $q \in \Delta(A)$ be given. The set of approachable targets is $$\mathcal T(q)=\{y\mid\,y=\sum_{l,k\in A} p_l q_k M_{lk}, \forall p\in \Delta(A)\}.$$
Furthermore, there exists a partitioning $R_1, \ldots, R_{n}$  such that the approachable strategies are markovian and bang-bang:
\begin{equation}\sigma(x)=\left\{\begin{array}{ll}
a_i=k & \mbox{if $x \in R_k:=\{\xi| \, (\xi-y)^T (u(k,q)-y) \leq 0\}$}\\ 
a_i\not = k & \mbox{otherwise.}\end{array}\right.\end{equation} 
\end{theorem}

\proof  Sketch. 
(sufficiency) Let  $y\in \mathcal T(q)$. Rewrite as $y= \sum_{l,k\in A} p_l q_k M_{lk}$ where where $p,q \in \Delta(A)$. Let us also take $\Phi = \{y(t)\}$. 

 
Then for every $x \in X$, taking $\lambda = \frac{x-y}{\|x-y\|} \in \mathbb R^m$ the value of the projected game satisfies
\begin{equation}\label{apprm1} 
\left\{\begin{array}{lll}
[x(t) - y]^T [x(t)-y] + \|x-y\| \mathbb E val_x(\lambda)\leq 0, \quad \forall \ t. \\
\mathbb E val_x(\lambda) := 
 \min_{a_i(t)}  \lambda^T f\Big(u(a_i(t),q(t) ),x\Big)\\
\end{array}\right.
\end{equation}

(necessity) Let $y\not \in \mathcal T(q)$. Then the above does not hold.  
 {\bf Q.E.D.}

\bigskip

In the problem at hand, one additional challenge is that $q$ must be \emph{self-confirmed}. This means that the mixed strategy $q$ entering the computation of the expected value of the projected games $\mathbb E val_x(\lambda)$ must reflect the current state distribution. In formulas, this corresponds to expanding (\ref{apprm1}) as follows:
\begin{equation}\label{apprm11} 
\left\{\begin{array}{lll}
[x(t) - y]^T [x(t)-y] + \|x- y\| \mathbb E val_x(\lambda)\leq 0, \quad \forall \ t. \\
\mathbb E val_x(\lambda) := 
 \min_{a_i(t)}  \lambda^T f\Big(u(a_i(t),q(t) ),x\Big)\\
q \in \Delta(A) \,  s.t. \, q_k = \int_{R_k} \rho(x,t) dx, \\ \qquad 
R_k:=\{\xi| \, (\xi-y)^T (u(k,q)-y) \leq 0\} \, \forall k \in A.
\end{array}\right.
\end{equation}
In the rest of the paper we look for self-confirmed solutions, which we call equilibria.
\subsection{The mean field game}\label{derive}
Let us denote by $v(x, t)$ the  value of the optimization problem starting
from time $t$ at state $x$. The first step is to show that the
problem results in the following mean field game system for
the unknown scalar functions $v(x, t)$, and $\rho(x, t)$ when each group behaves according to (\ref{cost}):

\begin{equation}
\label{eq:meanfieldcor}
\left\{
\begin{array}{l}
\displaystyle
\partial_t v(x,t)+\inf_{a_i } \left\{f(u(a_i,q),x) \partial_x v(x,t)+g(x,y)\right\} 
=0\ \mbox{in } \mathbb{R}^m\times[0,T[,\\
\\
\displaystyle
v(x,T)=\Psi(x,y)\ \forall\ x\in\mathbb{R}^m,\\
\displaystyle
\\
\partial_t \rho(x,t) + div(\rho(x,t) \cdot f(u(a_i^*,q),x)) 
= 0,\\
\\
\rho(0)= \rho_0,\\
\end{array}
\right.
\end{equation}
where $a_i^*(t,x)$ and $q$ are the optimal time-varying state-feedback
controls of players $i$ and $j$, respectively, obtained as
\begin{equation}
\label{eq:meanfield1new}
\left\{
\begin{array}{l}
a_i^* = \sigma(x) \in \arg \min_{a_i \in A_i}\{f(u(a_i,q),x) \partial_x v(x,t)+g(x,y)\},\\
\\
q \in \Delta(A) \,  s.t. \, q_k = \int_{R_k} \rho(x,t) dx, \\ \qquad  R_k:= \{x \in \mathbb R^m| \, \sigma(x)=k\}, \, \forall k \in A.\end{array}
\right.
\end{equation}

The mean field game system (\ref{eq:meanfieldcor}) appears in the form of two coupled PDEs intertwined in a forward-backward way. 
The first equation in (\ref{eq:meanfieldcor}) is the \emph{Hamilton-Jacobi-Bellman} (HJB) 
equation with variable $v(x,t)$ and parametrized in $\rho(\cdot)$. 
Given the boundary condition on final state (second equation in (\ref{eq:meanfieldcor})), and assuming a given population behavior captured by $\rho(\cdot)$,
the HJB equation is solved backwards and returns the value function and best-response behavior of the individuals (first equation in (\ref{eq:meanfield1new})) as well as 
the worst adversarial response (second equation in (\ref{eq:meanfield1new})). The HJB equation is coupled with a second PDE, known as \emph{Fokker-Planck-Kolmogorov (FPK)} (third equation in (\ref{eq:meanfieldcor})), defined on variable $\rho(\cdot)$ and parametrized in $v(x,t)$. Given the boundary condition on initial distribution $\rho(0)= \rho_0$ (fourth equation in (\ref{eq:meanfieldcor})), and assuming a given individual behavior described by $u^*$, the FPK equation is solved forward  and returns the population behavior time evolution $\rho(t)$.

Let condition (\ref{apprm}) hold true. Now, for given $x$, take for $\lambda$ the value $\lambda(\partial_x v)=\frac{\partial_x v(x,t)}{\|\partial_x v(x,t)\|}$ which is the gradient direction on $x$. Then, we can introduce the expected value of the \emph{projected anti-gradient game}  $$\mathbb E val_x[\partial_x v(x,t)]:=\lambda(\partial_x v)^T f(u_i(a_i^*,q),x).$$

We can then establish the following result. 
\begin{theorem}\label{polys} \textbf{(Self-confirmed equilibria)}
Let condition (\ref{apprm}) hold true. Then, the mean-field game formulation of Problem \ref{prob1}   is
\begin{equation}\label{cbl22}
\left\{\begin{array}{lll}
\partial_t v(x,t) +  \|\partial_x v \| \mathbb Eval_x[\partial_x v]
+ \frac{1}{2} (y(t) - x)^T Q (y(t) - x)
=0,\\
\; \ \mbox{ in } \mathbb{R}^m\times[0,T[,\\
\\
v(x,T)= \Psi (y(T),x),\ \mbox{in } \mathbb{R}^m,\\
\\ 
\partial_t \rho(x,t)  + div(\rho(x,t) \cdot f(u_i(a_i^*,q)) 
=0,
\ \mbox{in } \mathbb{R}^m\times[0,T[,\\
\\
\rho(x,0)= \rho_0(x)\ \mbox{in } \mathbb{R}^m.
\end{array}\right.
\end{equation}

Furthermore, the optimal controls for players 1 and 2  are
\begin{equation}\label{optul2}
\left\{\begin{array}{lll}
a_i^* = \sigma(x) \in   \arg \min_{a_i\in A_i} \lambda(\partial_x v)^T f(u(a_i,q),x) \\
q \in \Delta(A) \,  s.t. \, q_k = \int_{R_k} \rho(x,t) dx, \\ \qquad  
R_k:=\{\xi| \, (\xi-y)^T (u(k,q)-y) \leq 0\} \, \forall k \in A,\\
\sigma(x)=k,  \mbox{such that $x \in R_k$.} 
\end{array}\right.
\end{equation}

\end{theorem}

\proof{}  
Due to the bilinear structure of $f$, we can deduce that the best-response strategy $u^*$ and worst adversarial disturbance $w^*$ are 
on a vertex of the associated simplices in $\mathbb R^p$ and $\mathbb R^q$, respectively. This corresponds to saying
that both strategies are \emph{pure strategies}. We recall here that pure strategies are such that each player chooses as a result a  single predetermined action, in contrast with \emph{mixed strategies} where players select probabilities on actions and at time of play a random mechanism consistent with the selected probability distribution determines the actual action. A consequence of this is that  the mean field equilibrium, if exists, is  in pure strategies as well. 

We can   rewrite the value of the anti-gradient projected game as 
$$ \mathbb E val_x[ \partial_x v] = \inf_{l\in A}  \sum_{k \in A} \frac{1}{t} (q_k M_{lk} - x)^T  \lambda(\partial_x v),$$
Best responses and adversarial strategies are then given by
$$ a_i^* = \arg \min_{l \in A} \sum_{k \in A} \frac{1}{t} (q_k M_{lk} - x)^T  \lambda(\partial_x v).$$

With the above definition of $\mathbb E val_x[\partial_x v]$ in mind, the Hamilton-Jacobi part of (\ref{eq:meanfieldcor}) can be rewritten as
\begin{eqnarray}\nonumber
\partial_t v +  \|\partial_x v \| \mathbb E val_x [\partial_x v] +  \frac{1}{2}\left(y(t) - x(t)\right)^T Q \left(y(t) - x(t)\right) 
\label{eq:HJ0} 
=0\ \mbox{in } \mathbb{R}^m\times[0,T[, \\v(x,T)=\Psi(x)\ \forall\ x\in\mathbb{R}^m.
\end{eqnarray}

It is left to observe that $f(u^*,w^*)=A_{i^*j^*}$ and proves the third equation (FPK equation).
 {\bf Q.E.D.}

In principle, to find the optimal control input we need to solve
the two coupled PDEs in (\ref{cbl22})  in $v$ and $\rho$ with given boundary conditions (second and last conditions).

\subsection{Existence and nonuniqueness of equilibria}
In this section we investigate existence and nonuniqueness of equilibria. To do this, we analyze the time-dependence of an estimate error $\nu(t)$, which accounts for the deviation between an estimated  density $q(t)$ and a current one $\tilde q(t)$ at time $t$: 
$$\nu(t)= q(t) - \tilde q(t),$$
where
\begin{equation}\label{til}
\left\{\begin{array}{ll}
\tilde q_k(t) = \int_{R_k} \rho(x)dx \\
R_k:=\{\xi| \, (\xi-y(t))^T (u(k,q)-y(t)) \leq 0\}. 
\end{array}\right.\end{equation}
Observe that the time-dependence of $ \tilde q(t)$ enters in the above through the time-varying nature of the target $y(t)$.
Now, according to our procedure, we wish to hypothesize a pair $(p,q)$, which constitutes the input,  and obtain a new density $\tilde q(p,q)$ as an output. 
To see this, from $y=\sum_{l,k\in A} p_l q_k M_{lk}, \forall p,q\in \Delta(A)$
the expression (\ref{til}) can be rewritten as 
\begin{equation}
\left\{\begin{array}{ll}
\tilde q_k(p,q) = \int_{R_k} \rho(x)dx,\\
R_k:=\{\xi| \, (\xi-\sum_{l,k\in A} p_l q_k M_{lk})^T (u(k,q)-\sum_{l,k\in A} p_l q_k M_{lk}) \leq 0\}.
\end{array}\right.\end{equation}
Eventually, the procedure should return a fixed point. In other words,  if we think of an equilibrium as the pair $(p^*,q^*)$ such that $\nu(p^*,q^*) = 0$, existence of an equilibrium is now related to existence of a fixed point for the above procedure, i.e.,  $$\tilde q(p_1,q_1)= q.$$
The above means that, given a $(p,q)$ as input to our procedure, the output $\tilde q(p,q)$ coincides with the hypothesized density $q$. It is natural to represent the above algorithmic procedure, as a continuous-time dynamical system and thus to relate convergence to a fixed point to the asymptotic stability of the dynamics. The next assumption introduces conditions for the asymptotic stability to hold.  

\begin{assumption}\label{1}
There exists a pair $(\dot p,\dot q)$ such that 
\begin{equation}\label{cond}
\left[\begin{array}{c} 
- \partial_{p} \tilde q_1 \dot p + \dot q_1 - \partial_{q} \tilde q_1 \dot q  \\
\vdots \\
- \partial_{p} \tilde q_i \dot p + \dot q_i - \partial_{q} \tilde q_i \dot q
 \\ \vdots \\ 
  - \partial_{p} \tilde q_m \dot p + \dot q_m - \partial_{q} \tilde q_m \dot q
  \end{array} \right] := (- \partial_{p} \tilde q_i \dot p + \dot q_i - \partial_{q} \tilde q_i \dot q
)_{i=1,\ldots,m} \leq - \kappa (q - \tilde q).\end{equation}
\end{assumption}
The above describes the possibility of varying $(p,q)$ in order to reduce the estimate error $\nu$, whatever the current error is. The next result establishes the existence of an equilibrium based on the above condition. 

\begin{theorem}\label{thm:ex} {\bf(existence)}
Let Assumption \ref{1} hold. Then,  the estimate error decays exponentially fast, i.e.
$$\nu(t) \leq e^{-\kappa t}\nu(0).$$
\end{theorem}
\proof{ 
This proof is based on a Lyapunov stability approach. 
In particular, let us introduce a quadratic (in the error) Lyapunov function $$\mathcal L = \frac{1}{2}\nu^T \nu,$$ and show that its derivative is strictly negative. The time derivative can be decomposed as sum of two terms involving the gradient of $\mathcal L$ with respect to the two variables $p$ and $q$. More specifically, 
\begin{equation}
\begin{array}{ll}
\dot {\mathcal L} = (\partial_{p} \mathcal L)^T \dot p +  
(\partial_{q} \mathcal L)^T \dot q\\
=\nu^T \dot \nu = (q - \tilde q)^T \Big[ \left((\partial_{p} \nu_i)^T \dot p\right)_{i=1,\ldots,m} + \left((\partial_{q} \nu_i)^T \dot q \right)_{i=1,\ldots,m} \Big] \\
= (q - \tilde q)^T (- \partial_{p} \tilde q_i \dot p + \dot q_i - \partial_{q} \tilde q_i \dot q )_{i=1,\ldots,m}. 
\end{array}\end{equation} 
From condition (\ref{cond}), we also have $$\dot {\mathcal L} \leq - \kappa (q- \tilde q)^T(q- \tilde q) = - \kappa \nu^T \nu,$$
 which proves the thesis. {\bf Q.E.D.}

}
Essentially the above theorem shows that if we let the algorithm run for a long time the estimate error asymptotically converges to zero, namely,  $$\lim_{t \rightarrow \infty} \nu = 0,$$ which proves the existence of an equilibrium. 

We are now in the position to study nonuniqueness of equilibria. In particular, we provide a variational condition under which the equilibrium is nonunique.

\begin{theorem}\label{thm:uq} {\bf (nonuniqueness)}
Starting at an equilibrium where $\mathcal L=0$, if for all $\lambda \in \mathbb R^m$, $\|\lambda\|=1$ we have
\begin{equation}
\begin{array}{ccc}
\min_{\dot p,\dot q} \lambda^T \dot \nu =\min_{\dot p,\dot q} \lambda^T (- \partial_{p} \tilde q_i \dot p + \dot q_i - \partial_{q} \tilde q_i \dot q )_{i=1,\ldots,m}\\
< 0 <  \max_{\dot p,\dot q} \lambda^T \dot \nu = \max_{\dot p,\dot q} \lambda^T (- \partial_{p} \tilde q_i \dot p + \dot q_i - \partial_{q} \tilde q_i \dot q )_{i=1,\ldots,m},\end{array}\end{equation}
then there exists a $(\dot p,\dot q)$ such that $\dot{\mathcal L}=0$ and thus the current equilibrium is nonunique. 
\end{theorem}
\proof{There exists a  $(\dot p,\dot q)$ such that  $$\tilde q (p+\dot p dt,q+\dot q dt)= q+\dot q dt.$$ The above also means that the error $$\nu= \tilde q (p+\dot p dt,q +\dot q dt) - (q + \dot q dt )= 0.$$ {\bf Q.E.D.} }
\subsection{Solution of the mean field game}

This section investigates on the microscopic dynamics of every player given an equilibrium $(p,q)$ and the corresponding target which is common prior, where the target is denoted by 
$$y=\sum_{l,k\in A} p_l q_k M_{lk}.$$ 
As a result we obtain that such  a dynamics is a  ``potential'' one, in the sense that every player's current average payoff, which we can call \emph{state} of the player, describes a trajectory along the anti-gradient of a potential function, the latter being the value function of the mean-field game introduced earlier. To this purpose, let us denote by $e(t)$  the deviation between the target $y$ that every player wishes to approach, and the current average payoff $x(t)$, namely $$e(t) =y-x(t).$$  
Given that our running cost is quadratic, from dynamic programming, it is natural to assume that the value function has also a quadratic structure. This is a recurrent approach which needs an a posteriori verification of the consistency of the quadratic assumption. 
In particular, let us assume that the upper bound for the value function takes the form
\begin{equation}\label{phi}
\phi(x,t) = \frac{1}{2} e^T \Phi_t e,
\end{equation} 
where $\Phi_t$ is an opportune matrix which is positive definite, i.e.,  $\Phi_t >0$. 
Likewise, consider a quadratic function for the terminal penalty, namely, $$\Psi(x) = \frac{1}{2} e(T)^T \psi e(T).$$

Then, the HJB equation  in (\ref{eq:HJ0}) can be rewritten as
\begin{eqnarray}\nonumber
\partial_t \phi(x,t) +  \|\partial_x \phi(x,t) \| \mathbb E val_x [\partial_x \phi(x,t)] +  \frac{1}{2} e(t)^TQ e(t) 
\label{eq:HJ00} 
=0\ \mbox{in } \mathbb{R}^m\times[0,T[, \\ \phi(x,T)=\Psi(x)\ \forall\ x\in\mathbb{R}^m.
\end{eqnarray}
Substituting the expression (\ref{phi}) for the value function in (\ref{eq:HJ00}) we obtain
\begin{eqnarray}\nonumber
\frac{1}{2}e(t)^T \dot \Phi_t e(t) -\frac{1}{2} e(t)^T \Phi_t e(t) +  \frac{1}{2}e(t)^T Q e(t) 
\label{eq:HJ0} 
=0\ \mbox{in } \mathbb{R}^m\times[0,T[, \\  \frac{1}{2} e(T)^T \Phi_T e(T)=\Psi(x)\ \forall\ x\in\mathbb{R}^m.
\end{eqnarray}
The advantage of writing the HJB as above is in that all terms are explicitly written as quadratic terms in the error $e(t)$. Considering that the HJB has to hold true for every $e(t)$, we can drop $e(t)$ and thus we have an expression in the only matrix variable $\Phi_t$ as displayed next:
\begin{equation}\nonumber
\left\{\begin{array}{lll} 
\dot \Phi_t  - \Phi_t +   Q =0\ \mbox{in } [0,T[, \\ 
  \Phi_T= \psi \ \forall\ x\in\mathbb{R}^m.\end{array}\right.
\end{equation}
The above has the form of a classical differential Riccati equation which can be solved backwardly given the boundary conditions on the matrix in the terminal penalty, $\Phi_T= \psi$.
We can use such a result to analyze the microscopic dynamics of each player as detailed in the next subsection.

\subsubsection{Microscopic model}
Every single player is characterized by the following system of equations involving the evolution of the average payoff  (first equation), its best-response (second equation), and the expression for the density (third equation):
\begin{equation}\label{optul222}
\left\{\begin{array}{lll}
dx(t)= \frac{1}{t} \left(\sum_{k \in A} q_k M_{a^*k} - x(t) \right) dt,\\
a^*(x,t)=   \arg \min_{a\in A} (\Phi_t e(t))^T \left(\sum_{k \in A} q_k M_{ak} - x(t) \right), \\
q \in \Delta(A) \,  s.t. \, q_k = \int_{R_k} \rho(x,t) dx, \\ \qquad  R_k:= \{x \in \mathbb R^m| \, \sigma(x,t)=k\}, \, \forall k \in A.
\end{array}\right.
\end{equation}
Note that the expression for the best-response is obtained from (\ref{optul2}) where $\partial_x v$ is now replaced by $\Phi_t e(t)$. This is a straightforward consequence from assuming the value function quadratic as in (\ref{phi}).  

Let $t=e^s$ then

$$\dot x(s)= \sum_{k \in A} q_k M_{a^*k} - x(s) = u(a^*,q) - x(s).$$

For all $x$ the supporting hyperplane $H:=\{\xi| \, (\xi-y)^T (u(a^*,q)-y) = 0\}$ separates $x$ from 
$u(a^*,q)$, i.e., 

$$(x-y)^T (u(a^*,q) -y) = (x-y)^T (\sum_{k \in A} q_k M_{a^*k}-y) \leq 0.$$

Then from Theorem \ref{app_princ} approachability follows.

\section{Application: Regret and Bayesian equilibrium}\label{regret}
Perhaps the leading application of games with vector payoffs is in the study of regret-based dynamics, to which we now turn.

\subsection{Regret targeting in classical two-player games}
Given a symmetric normal-form game with common action set $A$ and symmetric payoff function $\pi:A\rightarrow\mathbb{R}$, let the \emph{regret} of player $i$ from not having played action $k\in A$ under action profile $\alpha\in A^2$ be
$$r(k,\alpha)=\pi(k,\alpha_{-i})-\pi(\alpha_i,\alpha_{-i}).$$
A straightforward way to justify the vector payoffs introduced earlier is to make them coincide with the regret vector associated to each action profile, i.e.
$$u(\alpha):=\Big(r(k,\alpha)\Big)_{k \in A}.$$
In Hart and Mas-Colell \cite{HM03}, approachability of the nonpositive orthant implies convergence to Nash equilibrium under such payoffs. This is no longer true for $1$st-moment approachability, which drives \emph{expected}---rather than maximum---regret to zero, so that some deviations could still have positive regret.

In the following, we turn standard games like the Prisoners' Dilemma, coordination games and Hawk--Dove games into games with regret vectors of type
$$\begin{array}{cc}
\begin{tabular}{|c|c|c|}
  \hline
       & Left & Right \\
  \hline
  Top & $\left(\begin{array}{c} 0 \\ a \end{array}\right)$ & $\left(\begin{array}{c} 0 \\ b \end{array}\right)$ \\
  \hline
  Bottom & $\left(\begin{array}{c} -a \\ 0 \end{array}\right)$ & $\left(\begin{array}{c} -b \\ 0 \end{array}\right)$ \\
  \hline
\end{tabular}
\end{array}$$
and analyse the resulting dynamics of a population targeting expected regret.

\begin{example} \textbf{(Prisoners' Regret)}
Consider again the Prisoners' Dilemma, and the following bimatrix, which represents the regret vector of player 1:
$$\begin{array}{cc}
\begin{tabular}{|c|c|c|}
  \hline
       & Cooperate & Defect \\
  \hline
  Cooperate & $\left(\begin{array}{c} 0 \\ 1 \end{array}\right)$ & $\left(\begin{array}{c} 0 \\ 1 \end{array}\right)$ \\
  \hline
  Defect & $\left(\begin{array}{c} -1 \\ 0 \end{array}\right)$ & $\left(\begin{array}{c} -1 \\ 0 \end{array}\right)$\\
  \hline
\end{tabular} 
\end{array}$$
Putting ourselves in the position of the Row player, and supposing that the Column player is randomly extracted from the population, 
we have that if Column is playing $D$, then if Row switched from $D$(efect) to $C$(ooperate), he would lose his payoff of $1$, whereas if he stuck to $D$ the regret would be $0$. This explains the vector payoff $(-1,0)$ for the action profile $(D,D)$. Likewise, if Row switched from $C$ to $D$ he would earn a payoff of $1$, in comparison with a regret of $0$ when sticking to $C$. This is represented by the regret vector $(0,1)$ for the action profile $(C,D)$.  The reasoning would be analogous if Column were to play $C$. 
Note that at the pure Nash equilibrium $(D,D)$ the regret vector is component-wise nonpositive.
\begin{figure} [htb]
\centering
\def\svgwidth{0.4\columnwidth}
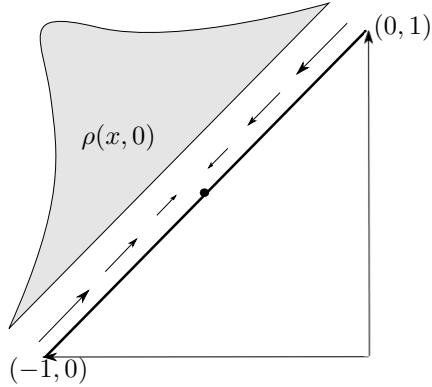
\caption{Regret space of Prisoners' dilemma: State space $X=conv\{(-1,0),(0,1)\}$ (solid line), initial distribution $\rho(x,0)$ (grey area), and vector field $dx(t)$ converging to $y=(-0.5,0.5)$.}  \label{fig:PD}
\end{figure}

 Figure~\ref{fig:PD} depicts the state space $X=conv\{(-1,0),(0,1)\}$ (solid line) in the case with an initial distribution $m(x,0)$ (grey area) of players. The arrows indicate the vector field $dx(t)$ if every player in state $x \in conv\{(-1,0),(-1/2,-1/2)\}$ cooperates, i.e.\ $a_i=1$ and every player in state $x \in conv\{(0,1),(-1/2,-1/2)\}$ defects. The vector field is such that eventually all players converge to the target$y=(-1/2,1/2)$. Consequently, the distribution converges asymptotically to a Dirac impulse in $y$.
\end{example}

\begin{example} \textbf{(Coordination game)}
Consider now the coordination game in the bimatrix on the left, with associated regret-vector game on the right:
$$\begin{array}{cc}
\begin{tabular}{|c|c|c|}
  \hline
       & Mozart & Mahler \\
  \hline
  Mozart & $(2,2)$ & $(0,0)$ \\
  \hline
  Mahler & $(0,0)$ & $(1,1)$ \\
  \hline
\end{tabular} $~~~~$ \begin{tabular}{|c|c|c|}
  \hline
       & Mozart & Mahler \\
  \hline
  Mozart & $\left(\begin{array}{c} 0 \\ -2 \end{array}\right)$ & $\left(\begin{array}{c} 0 \\ 1 \end{array}\right)$ \\
  \hline
  Mahler & $\left(\begin{array}{c} 2 \\ 0 \end{array}\right)$ & $\left(\begin{array}{c} -1 \\ 0 \end{array}\right)$\\
  \hline
\end{tabular} 
\end{array}$$
\begin{figure} [htb]
\centering
\def\svgwidth{0.5\columnwidth}
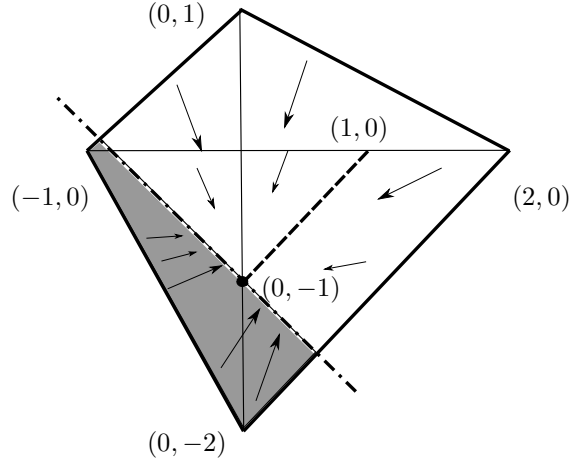
\caption{Regret space of the coordination game: State space $X=conv\{(-1,0),(0,1),(0,-2),(2,0)\}$ (boundary in solid line), and vector field $dx(t)$ converging to $(1,0)$ (grey area) and $(0,-1)$ (white area), approachable point is $y=(0,-1)$, set of approachable points is $conv\{(1,0),(0,-1)\}$ (dashed line) with mixed population strategy $q=(\frac{2}{3},\frac{1}{3})$.}  \label{fig:Coord}
\end{figure}
In Fig.\ \ref{fig:Coord} we illustrate the state space $X=conv\{(-1,0),(0,1),(0,-2),(2,0)\}$ (the boundary is in solid line). With a target $y=(0,-1)$, suppose we have a distribution on actions $q=(2/3,1/3)$, i.e.\ $2/3$ of the population plays Mozart, then $u(1,q)= (0,-1)$ and  $u(2,q)= (1,0)$ (here $k=2$ means playing Mahler). The set of approachable points with mixed population strategy $q=(2/3,1/3)$ is $conv\{(1,0),(0,-1)\}$ (dashed line), namely, any point in the convex hull of $u(1,q)= (0,-1)$ and  $u(2,q)= (1,0)$. The arrows indicate the vector field $dx(t)$ if every player in state $x \in R_2:=\{\xi| \, (\xi-y)^T (u(2,q)-y) \leq 0\}$ (grey area) plays Mahler, namely, $a_i=\sigma(x)=2$. On the other hand, every player in state $x \in R_1:=\{\xi| \, (\xi-y)^T (u(1,q)-y) \leq 0\}$ (white area) plays Mozart, namely, $a_i=\sigma(x)=1$.  Obviously we need that the integral of the distribution $m$ over $R_2$ is consistent with the initial assumption, which means $q_2=\int_{R_2} \rho(x,t) dx =1/3$. If this occurs, the vector field is such that eventually all players converge to $y=(0,-1)$. Consequently, the distribution converges to a Dirac impulse in $y$. 
\end{example}

\begin{example} (Hawk--Dove game)
We can likewise transform the Hawk--Dove (or chicken) game 
on the left into the corresponding regret-vector game on the right:
%
$$\begin{array}{cc}
\begin{tabular}{|c|c|c|}
  \hline
       & Hawk & Dove \\
  \hline
  Hawk & $\Big( -1,-1\Big)$ & (4,0) \\
  \hline
  Dove & (0,4) & $\Big( 2,2\Big)$ \\
  \hline
\end{tabular} $~~~~$ \begin{tabular}{|c|c|c|}
  \hline
       & Hawk & Dove \\
  \hline
  Hawk & $\left(\begin{array}{c} 0 \\ 1 \end{array}\right)$ & $\left(\begin{array}{c} 0 \\ -2 \end{array}\right)$ \\
  \hline
  Dove & $\left(\begin{array}{c} -1 \\ 0 \end{array}\right)$ & $\left(\begin{array}{c} 2 \\ 0 \end{array}\right)$\\
  \hline
\end{tabular} 
\end{array}$$
We have two pure Nash equilibria $(Dove, Hawk)$ and $(Hawk, Dove)$, whose corresponding regret vectors are nonpositive.\end{example}

\begin{figure} [htb]
\centering
\def\svgwidth{0.6\columnwidth}
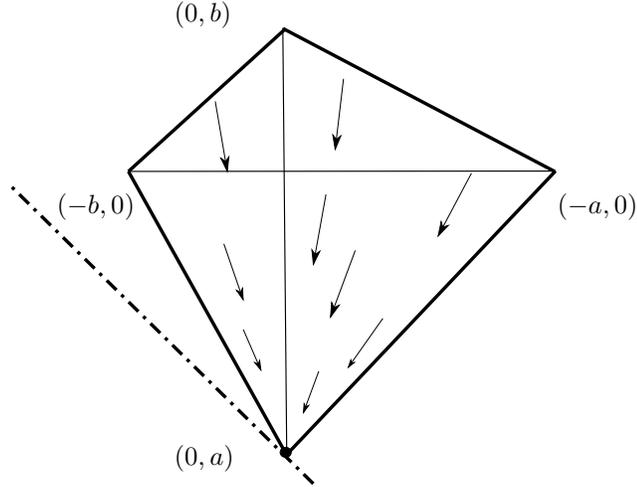
\caption{Regret space of parametric game with $a< 0 < b$: State space $X=conv\{(0,a),(-a,0),(-b,0),(0,b)\}$ (boundary in solid line), vector field $dx(t)$ converging to $(0,a)$ which is also an approachable vertex with population strategy $q=(1,0)$, supporting hyperplane $H$ (dot-dashed line) intersects $X$ only in one point (the vertex).}  \label{fig:ab1}
\end{figure}

More generally, let us now consider the  parametric game introduced earlier:
$$\begin{array}{cc}
\begin{tabular}{|c|c|c|}
  \hline
       & Left & Right \\
  \hline
  Top & $\left(\begin{array}{c} 0 \\ a \end{array}\right)$ & $\left(\begin{array}{c} 0 \\ b \end{array}\right)$ \\
  \hline
  Bottom & $\left(\begin{array}{c} -a \\ 0 \end{array}\right)$ & $\left(\begin{array}{c} -b \\ 0 \end{array}\right)$ \\
  \hline
\end{tabular}
\end{array}$$
Fig.\ \ref{fig:ab1} illustrates the state space $X=conv\{(0,a),(-a,0),(-b,0),(0,b)\}$ (the boundary is in solid line) where $a< 0 < b$. The target $y=(0,a)$ is in the negative orthant. 
 Here we consider a distribution on actions $q=(1,0)$, i.e.\ everybody plays $k=1$, then $u(1,q)= (0,a)$ and  $u(2,q)= (-a,0)$. The arrows indicate the vector field $dx(t)$ for which eventually all players converge to $y=(0,a)$. Consequently, the distribution converges to a Dirac impulse in $y$.
Note that the supporting hyperplane $H:=\{\xi| \, (\xi-y)^T (u(2,q)-y) = 0\}$ (dot-dashed line) intersects $X$ at only one point (the vertex), which is proven to be necessary for the vertex to be at the equilibrium. This will be explained in Theorem \ref{polys}. 

\begin{figure} [htb]
\centering
\def\svgwidth{0.4\columnwidth}
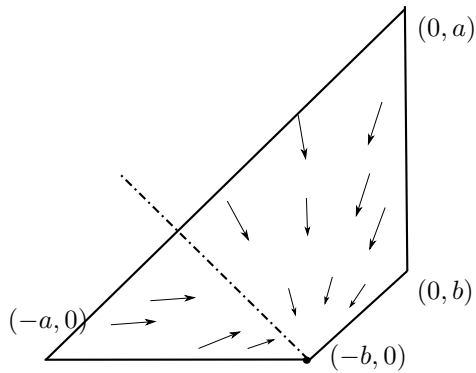
\caption{Regret space of parametric game with $0<b < a$: State space $X=conv\{(0,a),(-a,0),(-b,0),(0,b)\}$ (boundary in solid line), supporting hyperplane $H$ (dot-dashed line) passing through the vertex $(-b,0)$, vector field $dx(t)$ converging to $(0,b)$ left of $H$ and to $(-b,0)$ right of $H$, $conv\{(0,b),(-b,0)\}$ is set of approachable points with population strategy $q=(0,1)$, vertex $(-b,0)$ is not self-confirmed, while vertex $(0,a)$ is self-confirmed with population strategy $q=(1,0)$.}  \label{fig:ab2}
\end{figure} 

Fig.\ \ref{fig:ab2} depicts the state space $X=conv\{(0,a),(-a,0),(-b,0),(0,b)\}$ (the boundary is in solid line) where $0<b < a$. The target $y=(-b,0)$ is again in the negative orthant.  Here we consider a distribution on actions $q=(0,1)$, i.e.\ everybody plays $k=2$, then $u(1,q)= (0,b)$ and  $u(2,q)= (-b,0)$. The arrows indicate the vector field $dx(t)$ for which eventually all players converge to $y=(-b,0)$. Consequently, the distribution converges to a Dirac impulse in $y$.
However, there is an issue here related to the fact that the vertex $y$ is not at the equilibrium. To see this, note that the supporting hyperplane $H:=\{\xi| \, (\xi-y)^T (u(1,q)-y) = 0\}$ (dot-dashed line) partitions $X$ into two regions, which is proven to be necessary for the vertex not to be at the equilibrium. This will be explained in Theorem \ref{polys}.

\subsection{Maximum regret and Bayesian equilibrium}
Whilst $1$st-moment approachability gives interesting dynamics in population games based on regret then, it does not give convergence to Nash equilibrium. In this section, however, we show how the model can be applied to an incomplete-information setting to yield convergence to Bayesian equilibrium.

Suppose then that the continuous-time population game $\Gamma$ is based on a game of incomplete information; in particular, we are given a Harsanyi game $G$ (as described in \cite{Z12}) with state of the world $\omega=(s(\omega);t_1(\omega),t_2(\omega))$ chosen by Nature from a finite set $Y$ using a probability distribution $\theta$. Players then learn their own types $t_i(\omega)\in T_i$, choose actions $\beta_i$ from a common finite set $B(\omega)$, and receive symmetric payoffs $\varpi_i(\beta;\omega)$, $\beta=(\beta_1,\beta_2)$; the state of nature is $s(\omega)=(B(\omega),\varpi)$, $\varpi=(\varpi_1,\varpi_2)$. Each player $i$ then has a common finite set $\Sigma$ of ($T_i$-measurable) pure Bayesian strategies $\sigma_i:Y\rightarrow B(\omega)$, which we identify with the action set $A$ in our general framework.  Given a strategy profile $\sigma\in\Sigma^2$, let the vector payoffs be given by \emph{maximal regrets},
$$u(\sigma):=\Big(\max_{k \in\Sigma}r(k(\omega),\sigma(\omega))\Big)_{t_i\in T_i}.$$

Players are continuously rematched against new opponents to play this game $G$, and a new state of the world is chosen for each such matching; hence, each play of $G$ is one-shot in Nature, as distinct from repeated games of incomplete information (see \cite{AM95} and Ch.\ 14 of \cite{MSZ13}), where the opponents and state remain constant through time. $1$st-moment approachability of the nonpositive orthant in $\Gamma$ then implies that
$$\mathbb E_\theta\max_{k\in\Sigma}\varpi_i(k(\omega),\sigma_{-i}(\omega)) - \varpi_i(\sigma(\omega))\leq0.$$
But since the 
maximum of convex functions is convex, Jensen's inequality implies that the left-hand side is no less than
$$\max_{k\in\Sigma}\mathbb E_\theta\varpi_i(k(\omega),\sigma_{-i}(\omega)) - \mathbb E_\theta\varpi_i(\sigma(\omega)),$$
which is hence also nonpositive. Thus, we have a Nash equilibrium of the Harsanyi game, which is also a Bayesian equilibrium of the incomplete-information game by Harsanyi's \cite{H68} Theorem I.

For example, consider a game $G$ where each player's payoffs are randomly determined; with probability $1/2$, the Row player $R$ has the payoffs in the left-hand ``l'' matrix, and with probability $1/2$, she has the payoffs in the right-hand ``h'' matrix:
$$\begin{array}{cc}
\begin{tabular}{|c|c|c|}
  \hline
  $l$ & Opera & Football \\
  \hline
  Opera & $3$ & $1$ \\
  \hline
  Football & $0$ & $2$\\
  \hline
\end{tabular} $~~~~$ \begin{tabular}{|c|c|c|}
  \hline
  $h$ & Opera & Football \\
  \hline
  Opera & $1$ & $3$ \\
  \hline
  Football & $2$ & $0$\\
  \hline
\end{tabular} 
\end{array}$$
The Column player $C$'s payoffs are determined in a symmetric manner. Each player observes her own payoffs, but not those of her opponent. There are thus four possible states of the world $Y=\{\omega_{ll},\omega_{lh},\omega_{hl},\omega_{hh}\}$:
\begin{equation}
\left\{
\begin{array}{l}
\omega_{ll}=\left(s_{ll};[\frac{1}{2}\omega_{ll},\frac{1}{2}\omega_{lh}],[\frac{1}{2}\omega_{ll},\frac{1}{2}\omega_{hl}]\right) \\
\omega_{lh}=\left(s_{lh};[\frac{1}{2}\omega_{ll},\frac{1}{2}\omega_{lh}],[\frac{1}{2}\omega_{lh},\frac{1}{2}\omega_{hh}]\right) \\
\omega_{hl}=\left(s_{hl};[\frac{1}{2}\omega_{hl},\frac{1}{2}\omega_{hh}],[\frac{1}{2}\omega_{ll},\frac{1}{2}\omega_{hl}]\right) \\
\omega_{hh}=\left(s_{hh};[\frac{1}{2}\omega_{hl},\frac{1}{2}\omega_{hh}],[\frac{1}{2}\omega_{lh},\frac{1}{2}\omega_{hh}]\right),
\end{array}\right.
\end{equation}
each occurring with probability $1/4$. Furthermore, there are two possible types of each player,
$$\{R_l,R_h\}=\left\{\left[\frac{1}{2}\omega_{ll},\frac{1}{2}\omega_{lh}\right],\left[\frac{1}{2}\omega_{hl},\frac{1}{2}\omega_{hh}\right]\right\},$$
$$\{C_l,C_h\}=\left\{\left[\frac{1}{2}\omega_{ll},\frac{1}{2}\omega_{hl}\right],\left[\frac{1}{2}\omega_{ll},\frac{1}{2}\omega_{hl}\right]\right\},$$
and each player assigns probability $1/2$ to each of her opponents' possible types. Representing this situation as a Bayesian game, the Row player's vector payoffs are:
$$\begin{array}{cc}
\begin{tabular}{|c|c|c|c|c|}
  \hline
       & $O_l$, $O_h$ & $O_l$, $F_h$ & $F_l$, $O_h$ & $F_l$, $F_h$ \\
  \hline
  $O_l$, $O_h$ & $\left(\begin{array}{c} 3 \\ 1 \end{array}\right)$ & $\left(\begin{array}{c} 2 \\ 2 \end{array}\right)$ & $\left(\begin{array}{c} 2 \\ 2 \end{array}\right)$ & $\left(\begin{array}{c} 1 \\ 3 \end{array}\right)$ \\
  \hline
  $O_l$, $F_h$ & $\left(\begin{array}{c} 3 \\ 2 \end{array}\right)$ & $\left(\begin{array}{c} 2 \\ 1 \end{array}\right)$ & $\left(\begin{array}{c} 2 \\ 1 \end{array}\right)$ & $\left(\begin{array}{c} 1 \\ 0 \end{array}\right)$ \\
  \hline
  $F_l$, $O_h$ & $\left(\begin{array}{c} 0 \\ 1 \end{array}\right)$ & $\left(\begin{array}{c} 1 \\ 2 \end{array}\right)$ & $\left(\begin{array}{c} 1 \\ 2 \end{array}\right)$ & $\left(\begin{array}{c} 2 \\ 3 \end{array}\right)$ \\
  \hline
  $F_l$, $F_h$ & $\left(\begin{array}{c} 0 \\ 2 \end{array}\right)$ & $\left(\begin{array}{c} 1 \\ 1 \end{array}\right)$ & $\left(\begin{array}{c} 1 \\ 1 \end{array}\right)$ & $\left(\begin{array}{c} 2 \\ 0 \end{array}\right)$ \\
  \hline
\end{tabular}
\end{array}$$
where, for example, $O_l$, $F_h$ denotes the pure Bayesian strategy $\{\sigma_R(R_l)=\{\textrm{Opera}\},\sigma_R(R_h)=\{\textrm{Football}\}\}$. The Column player's payoffs are symmetric. This game has one pure-strategy equilibrium where Row plays $O_l$, $F_h$ and Column plays $O_l$, $O_h$, and a symmetric one where Row plays $O_l$, $O_h$ and Column plays $O_l$, $F_h$.

Now convert this game into one with maximal-regret payoffs:
$$\begin{array}{cc}
\begin{tabular}{|c|c|c|c|c|}
  \hline
       & $O_l$, $O_h$ & $O_l$, $F_h$ & $F_l$, $O_h$ & $F_l$, $F_h$ \\
  \hline
  $O_l$, $O_h$ & $\left(\begin{array}{c} 0 \\ 1 \end{array}\right)$ & $\left(\begin{array}{c} 0 \\ 0 \end{array}\right)$ & $\left(\begin{array}{c} 0 \\ 0 \end{array}\right)$ & $\left(\begin{array}{c} 1 \\ 0 \end{array}\right)$ \\
  \hline
  $O_l$, $F_h$ & $\left(\begin{array}{c} 0 \\ 0 \end{array}\right)$ & $\left(\begin{array}{c} 0 \\ 1 \end{array}\right)$ & $\left(\begin{array}{c} 0 \\ 1 \end{array}\right)$ & $\left(\begin{array}{c} 1 \\ 3 \end{array}\right)$ \\
  \hline
  $F_l$, $O_h$ & $\left(\begin{array}{c} 3 \\ 1 \end{array}\right)$ & $\left(\begin{array}{c} 1 \\ 0 \end{array}\right)$ & $\left(\begin{array}{c} 1 \\ 0 \end{array}\right)$ & $\left(\begin{array}{c} 0 \\ 0 \end{array}\right)$ \\
  \hline
  $F_l$, $F_h$ & $\left(\begin{array}{c} 3 \\ 0 \end{array}\right)$ & $\left(\begin{array}{c} 1 \\ 1 \end{array}\right)$ & $\left(\begin{array}{c} 1 \\ 1 \end{array}\right)$ & $\left(\begin{array}{c} 0 \\ 3 \end{array}\right)$ \\
  \hline
\end{tabular}
\end{array}$$
For instance, if Row is playing $F_l$, $O_h$ and Column is playing $O_l$, $O_h$, Row type $R_l$'s expected payoff is $0$, whereas he could have had $3$ by playing $O_l$, $O_h$, giving a maximal regret of $3$; similarly, type $R_h$'s payoff is $1$, whereas he could have had $2$ by playing $F_l$, $F_h$, giving a maximal regret of $1$. $1$st-moment approachability of the nonpositive orthant with these maximal-regret payoffs then implies Bayesian equilibrium. 

In this respect, from Theorem \ref{app_princ} we know that, for instance,  for any pure strategy $q$ we have
\begin{equation}
\mathcal T(q)=
\left\{\begin{array}{lll}
\{y\mid\,y \in conv( (0,1), (0,0), (3,1), (3,0))\}, & q=(1,0,0,0),\\
\{y\mid\,y \in conv( (0,0), (0,1), (1,0), (1,1))\}, & q=(0,1,0,0),\\
\{y\mid\,y \in conv( (0,0), (0,1), (1,0), (1,1))\}, & q=(0,0,1,0),\\
\{y\mid\,y \in conv( (1,0), (1,3), (0,0), (0,3))\}, & q=(0,0,0,1).
\end{array}\right.
\end{equation}
This means that for any pure strategy $q$ the origin $(0,0)$ is reachable and in particular the corresponding strategy is 
\begin{equation}\sigma(x)=\left\{\begin{array}{lll}
a_i=2 & \mbox{for all $x$}, & q=(1,0,0,0),\\ 
a_i=1 & \mbox{for all $x$}, & q=(0,1,0,0),\\ 
a_i=1 & \mbox{for all $x$}, & q=(0,0,1,0),\\ 
a_i=3 & \mbox{for all $x$}, & q=(0,0,0,1).\end{array}\right.\end{equation} 
However, note none of the above strategies corresponds to a self-confirmed equilibrium  according to Theorem \ref{polys}. Indeed, let us take for instance the first strategy, $a_i=2$, for all $x$ when $q=(1,0,0,0)$. But $a_i=2$, for all $x$ implies $R_2 = X$ and $R_1 = R_3 = R_4=   \emptyset$ which implies in turn $q=(0,1,0,0)$ and this contradicts the assumption $q=(1,0,0,0)$. We can repeat the same reasoning for any other strategy.

\section{Conclusion}\label{sec:conclusions}

This paper has combined approachability theory, evolutionary  games, and mean-field games in a unified framework. The game studied has a vector payoff, a large number of players, and admits classical mean-field game representation involving two coupled PDEs, the \emph{Hamilton-Jacobi-Bellman equation} and the \emph{advection equation}. We have highlighted multiple contributions. First, we coin the notion of \emph{1st-moment approachability} and analyze the corresponding convergence conditions. Second, we use the mean-field game to introduce the \emph{self-confirmed equilibrium}. Third we discuss on existence, non uniqueness, and stability of equilibria  as fixed points of the two PDEs. 

Future work involves the stochastic analysis of the same game in the presence of an additional Brownian motion in the dynamics. This would capture uncertainty or model-misspecification.  
In a different direction, we are interested in extending the study to the case where each player can adopt a mixed strategy, which would imply a new definition of density distribution on the space of mixed strategies; so far, the density distribution is defined on the space of pure strategies. A third development will be a further analysis of the connections with the Bayesian approach.



\bibliographystyle{plain}
\bibliography{References}


%
%

%

%
%
%


\end{document}

%% file: PDp.eps_tex
\begingroup%
  \makeatletter%
  \providecommand\color[2][]{%
    \errmessage{(Inkscape) Color is used for the text in Inkscape, but the package 'color.sty' is not loaded}%
    \renewcommand\color[2][]{}%
  }%
  \providecommand\transparent[1]{%
    \errmessage{(Inkscape) Transparency is used (non-zero) for the text in Inkscape, but the package 'transparent.sty' is not loaded}%
    \renewcommand\transparent[1]{}%
  }%
  \providecommand\rotatebox[2]{#2}%
  \ifx\svgwidth\undefined%
    \setlength{\unitlength}{393.85bp}%
    \ifx\svgscale\undefined%
      \relax%
    \else%
      \setlength{\unitlength}{\unitlength * \real{\svgscale}}%
    \fi%
  \else%
    \setlength{\unitlength}{\svgwidth}%
  \fi%
  \global\let\svgwidth\undefined%
  \global\let\svgscale\undefined%
  \makeatother%
  \begin{picture}(1,0.88453726)%
    \put(0,0){\includegraphics[width=\unitlength]{PDp.eps}}%
      \put(0.15,0.1){\color[rgb]{0,0,0}\makebox(0,0)[lb]{\smash{$\textstyle{(1,1)}$}}}%
 \put(1.0,0.04){\color[rgb]{0,0,0}\makebox(0,0)[lb]{\smash{$\textstyle{(4,0)}$}}}%
  \put(0.0,0.9){\color[rgb]{0,0,0}\makebox(0,0)[lb]{\smash{$\textstyle{(0,4)}$}}}%
 \put(0.8,0.75){\color[rgb]{0,0,0}\makebox(0,0)[lb]{\smash{$\textstyle{(3,3)}$}}}%
 \put(0.32,0.85){\color[rgb]{0,0,0}\makebox(0,0)[lb]{\smash{$\textstyle{(\frac{3}{2},\frac{7}{2})}$}}}%
  \put(0.55,0.01){\color[rgb]{0,0,0}\makebox(0,0)[lb]{\smash{$\textstyle{(\frac{5}{2},\frac{1}{2})}$}}}%
  \end{picture}%
\endgroup%

%% file: PD.eps_tex
\begingroup%
  \makeatletter%
  \providecommand\color[2][]{%
    \errmessage{(Inkscape) Color is used for the text in Inkscape, but the package 'color.sty' is not loaded}%
    \renewcommand\color[2][]{}%
  }%
  \providecommand\transparent[1]{%
    \errmessage{(Inkscape) Transparency is used (non-zero) for the text in Inkscape, but the package 'transparent.sty' is not loaded}%
    \renewcommand\transparent[1]{}%
  }%
  \providecommand\rotatebox[2]{#2}%
  \ifx\svgwidth\undefined%
    \setlength{\unitlength}{336.81464383bp}%
    \ifx\svgscale\undefined%
      \relax%
    \else%
      \setlength{\unitlength}{\unitlength * \real{\svgscale}}%
    \fi%
  \else%
    \setlength{\unitlength}{\svgwidth}%
  \fi%
  \global\let\svgwidth\undefined%
  \global\let\svgscale\undefined%
  \makeatother%
  \begin{picture}(1,0.98627789)%
    \put(0,0){\includegraphics[width=\unitlength]{PD.eps}}%
    \put(0.0,-0.04){\color[rgb]{0,0,0}\makebox(0,0)[lb]{\smash{$\textstyle{(-1,0)}$}}}%
    \put(1.0,0.9){\color[rgb]{0,0,0}\makebox(0,0)[lb]{\smash{$\textstyle{(0,1)}$}}}%
     \put(0.2,0.6){\color[rgb]{0,0,0}\makebox(0,0)[lb]{\smash{$\textstyle{\rho(x,0)}$}}}%
  \end{picture}%
\endgroup%

%% file: Coord.eps_tex
\begingroup%
  \makeatletter%
  \providecommand\color[2][]{%
    \errmessage{(Inkscape) Color is used for the text in Inkscape, but the package 'color.sty' is not loaded}%
    \renewcommand\color[2][]{}%
  }%
  \providecommand\transparent[1]{%
    \errmessage{(Inkscape) Transparency is used (non-zero) for the text in Inkscape, but the package 'transparent.sty' is not loaded}%
    \renewcommand\transparent[1]{}%
  }%
  \providecommand\rotatebox[2]{#2}%
  \ifx\svgwidth\undefined%
    \setlength{\unitlength}{324.321912bp}%
    \ifx\svgscale\undefined%
      \relax%
    \else%
      \setlength{\unitlength}{\unitlength * \real{\svgscale}}%
    \fi%
  \else%
    \setlength{\unitlength}{\svgwidth}%
  \fi%
  \global\let\svgwidth\undefined%
  \global\let\svgscale\undefined%
  \makeatother%
  \begin{picture}(1,0.93117359)%
    \put(0,0){\includegraphics[width=\unitlength]{Coord.eps}}%
    \put(0.2,-0.04){\color[rgb]{0,0,0}\makebox(0,0)[lb]{\smash{$\textstyle{(0,-2)}$}}}%
    \put(-0.1,0.5){\color[rgb]{0,0,0}\makebox(0,0)[lb]{\smash{$\textstyle{(-1,0)}$}}}%
    \put(1.0,0.5){\color[rgb]{0,0,0}\makebox(0,0)[lb]{\smash{$\textstyle{(2,0)}$}}}%
    \put(0.2,0.9){\color[rgb]{0,0,0}\makebox(0,0)[lb]{\smash{$\textstyle{(0,1)}$}}}%
    \put(0.6,0.65){\color[rgb]{0,0,0}\makebox(0,0)[lb]{\smash{$\textstyle{(1,0)}$}}}%
    \put(0.45,0.3){\color[rgb]{0,0,0}\makebox(0,0)[lb]{\smash{$\textstyle{(0,-1)}$}}}%
 \end{picture}%
\endgroup%

%% file: ab1.eps_tex
\begingroup%
  \makeatletter%
  \providecommand\color[2][]{%
    \errmessage{(Inkscape) Color is used for the text in Inkscape, but the package 'color.sty' is not loaded}%
    \renewcommand\color[2][]{}%
  }%
  \providecommand\transparent[1]{%
    \errmessage{(Inkscape) Transparency is used (non-zero) for the text in Inkscape, but the package 'transparent.sty' is not loaded}%
    \renewcommand\transparent[1]{}%
  }%
  \providecommand\rotatebox[2]{#2}%
  \ifx\svgwidth\undefined%
    \setlength{\unitlength}{385.309872bp}%
    \ifx\svgscale\undefined%
      \relax%
    \else%
      \setlength{\unitlength}{\unitlength * \real{\svgscale}}%
    \fi%
  \else%
    \setlength{\unitlength}{\svgwidth}%
  \fi%
  \global\let\svgwidth\undefined%
  \global\let\svgscale\undefined%
  \makeatother%
  \begin{picture}(1,0.83679403)%
    \put(0,0){\includegraphics[width=\unitlength]{ab1.eps}}%
    \put(0.3,0.04){\color[rgb]{0,0,0}\makebox(0,0)[lb]{\smash{$\textstyle{(0,a)}$}}}%
    \put(0.085,0.5){\color[rgb]{0,0,0}\makebox(0,0)[lb]{\smash{$\textstyle{(-b,0)}$}}}%
    \put(1.0,0.5){\color[rgb]{0,0,0}\makebox(0,0)[lb]{\smash{$\textstyle{(-a,0)}$}}}%
    \put(0.3,0.85){\color[rgb]{0,0,0}\makebox(0,0)[lb]{\smash{$\textstyle{(0,b)}$}}}%
   \end{picture}%
\endgroup%

%% file: ab2.eps_tex
\begingroup%
  \makeatletter%
  \providecommand\color[2][]{%
    \errmessage{(Inkscape) Color is used for the text in Inkscape, but the package 'color.sty' is not loaded}%
    \renewcommand\color[2][]{}%
  }%
  \providecommand\transparent[1]{%
    \errmessage{(Inkscape) Transparency is used (non-zero) for the text in Inkscape, but the package 'transparent.sty' is not loaded}%
    \renewcommand\transparent[1]{}%
  }%
  \providecommand\rotatebox[2]{#2}%
  \ifx\svgwidth\undefined%
    \setlength{\unitlength}{412.075bp}%
    \ifx\svgscale\undefined%
      \relax%
    \else%
      \setlength{\unitlength}{\unitlength * \real{\svgscale}}%
    \fi%
  \else%
    \setlength{\unitlength}{\svgwidth}%
  \fi%
  \global\let\svgwidth\undefined%
  \global\let\svgscale\undefined%
  \makeatother%
  \begin{picture}(1,0.98341433)%
    \put(0,0){\includegraphics[width=\unitlength]{ab2.eps}}%
    \put(1.02,0.18){\color[rgb]{0,0,0}\makebox(0,0)[lb]{\smash{$\textstyle{(0,b)}$}}}%
    \put(1.02,0.9){\color[rgb]{0,0,0}\makebox(0,0)[lb]{\smash{$\textstyle{(0,a)}$}}}%
    \put(-0.1,0.1){\color[rgb]{0,0,0}\makebox(0,0)[lb]{\smash{$\textstyle{(-a,0)}$}}}%
    \put(0.78,-0.0){\color[rgb]{0,0,0}\makebox(0,0)[lb]{\smash{$\textstyle{(-b,0)}$}}}%
  \end{picture}%
\endgroup%